\documentclass{amsart}

\usepackage{amssymb,amsthm,amsmath}
\usepackage{pstricks,pst-node,hyperref}

\begin{document}

\title[A combinatorial approach to the power of 2 on involutions]{A combinatorial approach to the power of 2 in the number of
involutions}
\author{Dongsu Kim}
\author{Jang Soo Kim}
\email[Dongsu Kim]{dongsu.kim@kaist.ac.kr}
\email[Jang Soo Kim]{jskim@kaist.ac.kr}
\address{Department of Mathematical Sciences\\
Korea Advanced Institute of Science and Technology\\
Daejeon 305-701, Korea}%
\thanks{The first author was supported by Basic Science Research Program through
the National Research Foundation of Korea(NRF) funded by the Ministry of
Education, Science and Technology (2009-0063183). The second author was supported by the second stage of the Brain Korea 21 Project, The Development Project of Human Resources in Mathematics, KAIST in 2009.}

\begin{abstract}
We provide a combinatorial approach to the largest power of $p$ in the
number of permutations $\pi$ with $\pi^p=1$, for a fixed prime number
$p$. With this approach, we find the largest power of $2$ in the
number of involutions, in the signed sum of involutions and in the
numbers of even or odd involutions.
\end{abstract}

\maketitle

% ----------------------------------------------------------------

\renewcommand\sectionautorefname{Section}
\newtheorem{thm}{Theorem}[section]
\newcommand\thmautorefname{Theorem}
\newtheorem{lem}{Lemma}[section]
\makeatletter
\renewcommand*{\c@lem}{\c@thm}
\makeatother
\newcommand\lemautorefname{Lemma}
\newtheorem{prop}{Proposition}[section]
\makeatletter
\renewcommand*{\c@prop}{\c@thm}
\makeatother
\newcommand\propautorefname{Proposition}
\newtheorem{cor}{Corollary}[section]
\makeatletter
\renewcommand*{\c@cor}{\c@thm}
\makeatother
\newcommand\corautorefname{Corollary}
\theoremstyle{definition}
\newtheorem{example}{Example}[section]
\makeatletter
\renewcommand*{\c@example}{\c@thm}
\makeatother
\newcommand\exampleautorefname{Example}
\newtheorem{defn}{Definition}[section]
\makeatletter
\renewcommand*{\c@defn}{\c@thm}
\makeatother
\newcommand\defnautorefname{Definition}
\newtheorem{conj}{Conjecture}[section]
\makeatletter
\renewcommand*{\c@conj}{\c@thm}
\makeatother
\newcommand\conjautorefname{Conjecture}
\newtheorem{question}{Question}[section]
\makeatletter
\renewcommand*{\c@question}{\c@thm}
\makeatother
\newcommand\questionautorefname{Question}
\theoremstyle{remark}
\newtheorem{remark}{Remark}[section]
\makeatletter
\renewcommand*{\c@remark}{\c@thm}
\makeatother
\newcommand\remarkautorefname{Remark}
\newtheorem*{note}{Note}

\newcommand\akn[3]{(#1;#2)_{#3}}

\def\iff.{if and only if}
\newcommand{\floor}[2]{\left\lfloor\frac {#1}{#2}\right\rfloor}
\newcommand{\fl}[1]{\left\lfloor {#1}\right\rfloor}
\newcommand{\ceil}[2]{\left\lceil\frac {#1}{#2}\right\rceil}
\def\S{\mathfrak{S}}
\def\I{\mathfrak{I}}
\def\N{\mathbb{N}}
\def\G{\mathfrak{G}}
\def\B{\mathfrak{B}}
\def\A{\mathcal{A}}
\def\seq#1#2#3{#1_1 #3 #1_2 #3\def\temp{#3}\def\comma{,}
 \if\comma\temp\ldots\else\cdots\fi #3 #1_{#2} }
\def\ao{\sigma_1}
\def\at{\sigma_2}
\def\ai{\sigma_i}
\def\wt{{\rm wt}}
\def\s{{\rm sign}}
\def\ord{{\rm ord}}
\def\inv{t}
\def\odd{\beta}

\newdimen\psunit\psunit=12pt
\newcount\a\newcount\b\newcount\c\newcount\d
\def\edgecolor{blue}
\def\ledge#1#2{\ncline[linecolor=\edgecolor]{#1}{#2}}
\psset{unit=\psunit}
\def\vput(#1)#2{\cnode(#1){8pt}{#2}\rput(#1){$v_{#2}$}}
\def\vbput(#1)#2{\cnode(#1){5pt}{#2}\rput(#1){$v_{#2}$}}
\def\vrput(#1)#2{\cnode[linecolor=red](#1){8pt}{#2}\rput(#1){$v_{#2}$}}
\def\iput(#1)#2{\cnode[linewidth=0.7pt](#1){5.5pt}{#2}\rput(#1){{\tiny $#2$}}}
\def\pput(#1)#2{\cnode[linewidth=0.7pt](#1){3.5pt}{#2}}
\def\spput(#1)#2{\cnode[linewidth=0.7pt](#1){1.7pt}{#2}}
\def\scpput(#1)#2{\cnode[fillstyle=solid,fillcolor=black,linewidth=0.7pt](#1){1.7pt}{#2}}
\def\eput(#1)#2{\cnode(#1){5pt}{#2}}
\def\uedge#1#2{\nccurve[linecolor=\edgecolor,angleA=60,angleB=120]{#1}{#2}}
\def\edge#1#2{\nccurve[linecolor=\edgecolor,angleA=300,angleB=240]{#1}{#2}}
\def\vb(#1,#2){\psset{unit=4pt}
\a=#1\b=#2\c=#1\d=#2
\multiply\a by3\multiply\b by3\multiply\c by3\multiply\d by3
\advance\c by3
\advance\a by-2\advance\b by-2\advance\c by2\advance\d by2
\psframe[framearc=1,linecolor=green!80!black](\a,\b)(\c,\d)
\psset{unit=13pt}}

\def\vbr(#1,#2){\psset{unit=4pt}
\a=#1\b=#2\c=#1\d=#2
\multiply\a by3\multiply\b by3\multiply\c by3\multiply\d by3
\advance\c by3
\advance\a by-2\advance\b by-2\advance\c by2\advance\d by2
\psframe[framearc=1,linecolor=red](\a,\b)(\c,\d)
\psset{unit=13pt}}
\def\xput(#1)#2{\cnode(#1){8pt}{#2}\rput(#1){$#2$}}
\def\labelput(#1)#2#3{\cnode[linewidth=0.7pt](#1){5.5pt}{#2}\rput(#1){{\tiny $#3$}}}

\section{Introduction}
The largest power of a prime in some well-known numbers has been studied
in many papers, for instance, see
\cite{Alter1973,Chowla1951,Deutsch2006,Goetgheluck1987,Grady1994,
  Ishihara2001,Knuth1989,Konvalinka2007,Ochiai1999,Postnikov2007}. In
this paper we are interested in the largest power of a prime in the
numbers of permutations with some conditions.

Let $\S_n$ denote the set of permutations of $[n]=\{1,2,\ldots,n\}$.
Let $p$ be a prime number and $n$ a positive integer. Let $\tau_p(n)$
denote the number of permutations $\pi\in\S_n$ such that $\pi^p=1$,
and let $\ord_p(n)$ denote the largest integer $k$ such that $p^k$
divides $n$.

In 1951, using recurrence relation with induction, Chowla, Herstein and Moore \cite{Chowla1951} proved that
$$\ord_2(\tau_2(n))\geq\floor n2-\floor n4.$$
Using generating function, Grady and Newman \cite{Grady1994} obtained, for any prime $p$,
\begin{equation}\label{grady}
\ord_p(\tau_p(n))\geq\floor np-\floor n{p^2}.
\end{equation}
Using $p$-adic analysis, Ochiai \cite{Ochiai1999} found the exact
value of $\ord_p(\tau_p(n))$ for prime numbers $p\leq 23$. Let $\inv_n$
denote $\tau_2(n)$, the number of involutions in $\S_n$. Ochiai's result
gives
\begin{equation}\label{ochiai}
\ord_2(t_n)=\floor n2-2\floor n4+\floor{n+1}4.
\end{equation}

In addition, Chowla at~el.~\cite{Chowla1951} considered the
sequence $\{t_n\mod m\}_{n\geq0}$ for a fixed integer $m$ and proved
that $m$ is a period of the sequence if $m$ is odd. We will prove that
in fact, it is the smallest period. If $m$ is even, then the
sequence is not periodic because $\inv_0=1$ but $\inv_n$ is even for
all $n\geq2$. However there is an integer $N$ such that $\{t_n\mod
m\}_{n\geq N}$ is periodic.

Our main results are in \autoref{sec:comb} and \autoref{tn}, where we
prove \eqref{grady} and \eqref{ochiai} using combinatorial arguments.
The weighted sum of involutions is considered in \autoref{sec:wt}.
In \autoref{even} we find $\ord_2$ of the signed sum of
involutions, the number of odd involutions, and the number of even
involutions. In
\autoref{period} we find the smallest $N$ such that $\{t_n\mod
m\}_{n\geq N}$ is periodic and find the smallest period of the
sequence when $m$ is even. We also consider the odd factor of the
number of involutions and prove that the smallest period of the
sequence $\{\inv_n/2^{\ord_2(\inv_n)}\mod 2^s\}_{n\geq0}$ is $2^{s+1}$
if $s\geq3$.

\section{A combinatorial proof}\label{sec:comb}
Let $\S_{n,p}$ denote the set of permutations $\pi\in\S_n$ with $\pi^p
=1$. For instance, for $p=2$ it is the set of all involutions in $\S_n$.
Each permutation in $\S_{n,p}$ is a product of disjoint $p$-cycles and
$1$-cycles. For example, for $\pi=38725614\in\S_{8,3}$, the disjoint
product is $(1,3,7)(2,8,4)(5)(6)$. A cycle usually consists of distinct 
integers, but we allow cycles to have repeated entries for convenience.

We define a {\em label map} $f_p:\{1,2,\dots,n\}\rightarrow\{1,2,\dots,\lfloor(n-1)/p\rfloor+1\}$ by $f_p(i)=
\lfloor(i-1)/p\rfloor+1$, extend it to cycles $\sigma=(s_1,\dots,s_j)$ by
$f_p(\sigma)=(f_p(s_1),\dots,f_p(s_j))$ which is regarded as a cycle with
repeated entries, and to $\S_{n,p}$ by
$$
f_p(\pi)=\{f_p(\sigma_1),\dots,f_p(\sigma_k)\}
$$
for $\pi=\sigma_1\sigma_2\cdots\sigma_k$ in the disjoint cycle notation.
Note that $f_p(\pi)$ is regarded as a multiset.

As a map defined on $\S_{n,p}$, $f_p$ induces an equivalence relation
$\sim$ on $\S_{n,p}$, namely $\pi\sim\tau$
if and only if $f_p(\pi)=f_p(\tau)$.

\begin{figure}
\begin{pspicture}(-3.46,-4)(3.46,3)
\iput(-2.6,2.5){1}\iput(-3.46,1){2}\iput(-1.73,1){3}
\iput(2.6,2.5){4}\iput(1.73,1){5}\iput(3.46,1){6}
\iput(0,-2){7}\iput(-0.87,-3.5){8}\iput(0.87,-3.5){9}
\ncline{->}{3}{5}
\ncline{->}{5}{7}
\ncline{->}{7}{3}
\nccurve[angleA=0,angleB=150]{->}{1}{6}
\nccurve[angleA=240,angleB=30]{->}{6}{8}
\nccurve[angleA=120,angleB=270]{->}{8}{1}
\nccurve[angleA=30,angleB=180,border=1pt]{->}{2}{4}
\nccurve[angleA=270,angleB=60,border=1pt]{->}{4}{9}
\nccurve[angleA=150,angleB=300,border=1pt]{->}{9}{2}
\end{pspicture}
\hspace{0.7cm}
\begin{pspicture}(-3.46,-4)(3.46,3)
\iput(-2.6,2.5){10}\iput(-3.46,1){11}\iput(-1.73,1){12}
\iput(2.6,2.5){13}\iput(1.73,1){14}\iput(3.46,1){15}
\iput(0,-2){16}\iput(-0.87,-3.5){17}\iput(0.87,-3.5){18}
\ncline{->}{14}{12}
\ncline{->}{16}{14}
\ncline{->}{12}{16}
\nccurve[angleA=0,angleB=150]{->}{10}{15}
\nccurve[angleA=240,angleB=30]{->}{15}{17}
\nccurve[angleA=120,angleB=270]{->}{17}{10}
\nccircle{->}{11}{.18cm}
\nccircle{->}{13}{.18cm}
\nccircle{->}{18}{.18cm}
\end{pspicture}
\hspace{0.7cm}
\begin{pspicture}(-3.46,-4)(3.46,3)
\iput(-2.6,2.5){19}\iput(-3.46,1){20}\iput(-1.73,1){21}
\iput(2.6,2.5){22}\iput(1.73,1){23}\iput(3.46,1){24}
\iput(-2.6,-2){25}\iput(-3.46,-3.5){26}\iput(-1.73,-3.5){27}
\iput(2.6,-2){28}\iput(1.73,-3.5){29}
\ncline{->}{19}{20}
\ncline{->}{20}{21}
\ncline{->}{21}{19}
\ncline{->}{28}{29}
\ncline{->}{29}{27}
\ncline{->}{27}{28}
\nccircle{->}{25}{.18cm}
\nccircle{->}{26}{.18cm}
\nccircle{->}{22}{.18cm}
\nccircle{->}{23}{.18cm}
\nccircle{->}{24}{.18cm}
\end{pspicture}
{$$\Downarrow f_3$$}
\begin{pspicture}(-3.46,-5)(3.46,4)
\labelput(-2.6,2.5){1}1\labelput(-3.46,1){2}1\labelput(-1.73,1){3}1
\labelput(2.6,2.5){4}2\labelput(1.73,1){5}2\labelput(3.46,1){6}2
\labelput(0,-2){7}3\labelput(-0.87,-3.5){8}3\labelput(0.87,-3.5){9}3
\ncline{->}{3}{5}
\ncline{->}{5}{7}
\ncline{->}{7}{3}
\nccurve[angleA=0,angleB=150]{->}{1}{6}
\nccurve[angleA=240,angleB=30]{->}{6}{8}
\nccurve[angleA=120,angleB=270]{->}{8}{1}
\nccurve[angleA=30,angleB=180,border=1pt]{->}{2}{4}
\nccurve[angleA=270,angleB=60,border=1pt]{->}{4}{9}
\nccurve[angleA=150,angleB=300,border=1pt]{->}{9}{2}
\end{pspicture}
\hspace{0.7cm}
\begin{pspicture}(-3.46,-5)(3.46,4)
\labelput(-2.6,2.5){10}4\labelput(-3.46,1){11}4\labelput(-1.73,1){12}4
\labelput(2.6,2.5){13}5\labelput(1.73,1){14}5\labelput(3.46,1){15}5
\labelput(0,-2){16}6\labelput(-0.87,-3.5){17}6\labelput(0.87,-3.5){18}6
\ncline{->}{14}{12}
\ncline{->}{16}{14}
\ncline{->}{12}{16}
\nccurve[angleA=0,angleB=150]{->}{10}{15}
\nccurve[angleA=240,angleB=30]{->}{15}{17}
\nccurve[angleA=120,angleB=270]{->}{17}{10}
\nccircle{->}{11}{.18cm}
\nccircle{->}{13}{.18cm}
\nccircle{->}{18}{.18cm}
\end{pspicture}
\hspace{0.7cm}
\begin{pspicture}(-3.46,-5)(3.46,4)
\labelput(-2.6,2.5){19}7\labelput(-3.46,1){20}7\labelput(-1.73,1){21}7
\labelput(2.6,2.5){22}8\labelput(1.73,1){23}8\labelput(3.46,1){24}8
\labelput(-2.6,-2){25}9\labelput(-3.46,-3.5){26}9\labelput(-1.73,-3.5){27}9
\labelput(2.6,-2){28}{10}\labelput(1.73,-3.5){29}{10}
\ncline{->}{28}{29}
\ncline{->}{29}{27}
\ncline{->}{27}{28}
\nccircle{->}{25}{.18cm}
\nccircle{->}{26}{.18cm}
\nccircle{->}{22}{.18cm}
\nccircle{->}{23}{.18cm}
\nccircle{->}{24}{.18cm}
\ncline{->}{19}{20}
\ncline{->}{20}{21}
\ncline{->}{21}{19}
\end{pspicture}
\caption{Visualization of $\pi$ and $f_3(\pi)$ in \autoref{ex}.}\label{fig:map}
\end{figure}
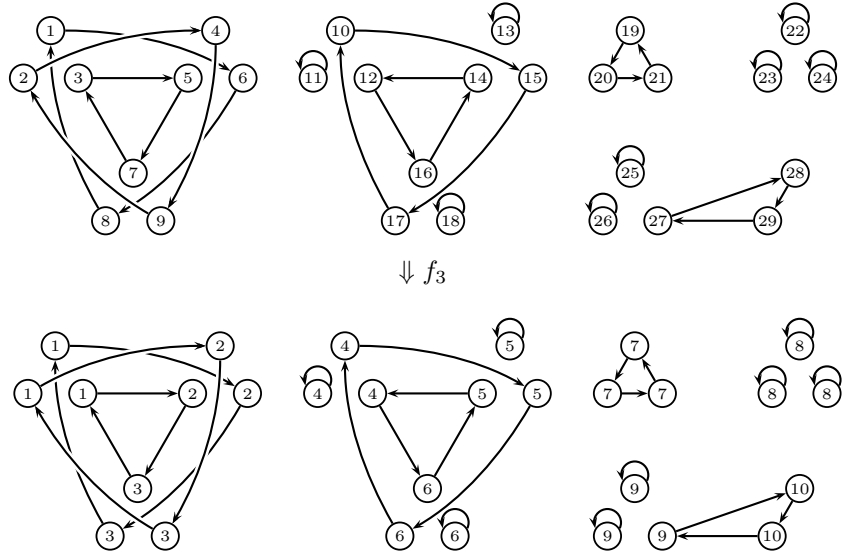

Fix a prime $p$, and let $n=pt+r$ with $0\leq r<p$.
A $p$-cycle $\sigma=(s_1,s_2,\dots,s_p)$ in some $f_p(\pi)$ is said to be
of \emph{type} $A$ if $s_1=s_2=\cdots=s_p$; of \emph{type} $B$ otherwise.
We are interested in the size of each equivalence class of $\sim$ on
$\S_{n,p}$. As a matter of fact, we need the size of some collections
of equivalence classes. An equivalence class may be represented as
a multiset of cycles with repeated entries from $\{1,2,\dots,t+1\}$.
In fact there are three kinds of cycles in the representation of 
equivalence classes:
$p$-cycles of type $A$, $p$-cycles of type $B$, and $1$-cycles.
A typical equivalence class is of the form
$\{A_1,\dots,A_i;B_1^{d_1},\dots,B_j^{d_j};C_1^{e_1},\dots,C_k^{e_k}\}$,
as a multiset, where $A$'s denote $p$-cycles of type $A$, $B$'s denote
those of type $B$ and $C$'s are $1$-cycles.
Since the multiplicities $e_1,\dots,e_k$ play a critical role,
we refine the form to
$\{A_1,\dots,A_i;B_1^{d_1},\dots,B_j^{d_j};C_1^{e_1},\dots,C_k^{e_k};D_1^p,
\dots,D_\ell^p\}$ with $e_1,\dots,e_k<p$, where $A$'s, $B$'s, $C$'s
are the same as before, while $D$'s are $1$-cycles.
We collect all equivalence classes
$\{A_1,\dots,A_i;B_1^{d_1},\dots,B_j^{d_j};C_1^{e_1},\dots,C_k^{e_k};D_1^p,
\dots,D_\ell^p\}$ with fixed $B$'s, $C$'s, and a fixed set of integers
appearing in either $A$'s or $D$'s. Let
$$
\{s_1,s_2,\dots,s_h;B_1^{d_1},\dots,B_j^{d_j};C_1^{e_1},\dots,C_k^{e_k}\}
$$
denote such a collection. The collection may be represented as
$$
\{s_1,s_2,\dots,s_h;E_1^{m_1},\dots,E_\ell^{m_\ell}\},
$$
where $E$'s denote either a $p$-cycle of type $B$ of multiplicity at most
$p$ or a $1$-cycle with multiplicity less than $p$.
Note that $\{s_1,\dots,s_h\}\subset[t]$ and each $i\in[t]
\setminus\{s_1,\dots,s_h\}$ appears exactly $p$ times in the collection
and $t+1$ appears exactly $r$ times. The distinct collections produce
a partition of $\S_{n,p}$, which in turn defines an equivalence relation,
denoted by $\sim'$. Let $\tilde f_p$ denote the quotient 
map corresponding to this equivalence relation.
\begin{example}\label{ex}
Let $\pi\in\S_{29,3}$ be the following permutation in cycle notation:
\begin{align*}
\pi=&(1,6,8)(2,4,9)(3,5,7)(10,15,17)(11)(12,16,14)(13)(18)\\
&(19,20,21)(22)(23)(24)(25)(26)(27,28,29).
\end{align*}
Then $f_3(\pi)=\{(1,2,3)^3,(4,5,6),(4),(4,6,5),(5),(6),(7,7,7),(8)^3,
(9)^2,(9,10,10)\}$.
The permutation $\pi$ belongs to an equivalence class
$$
\{(7,7,7);(1,2,3)^3,(4,5,6),(4,6,5),(9,10,10);(4),(5),(6),(9)^2;
(8)^3\}
$$
of the form $\{A_1,\dots,A_i;B_1^{d_1},\dots,B_j^{d_j};C_1^{e_1},\dots,
C_k^{e_k};D_1^p,
\dots,D_\ell^p\}$, which is a member of the collection
$$
\tilde f_3(\pi)=\{7,8;(1,2,3)^3,(4,5,6),(4,6,5),(9,10,10),(4),(5),(6),(9)^2\}
$$
We visualize this example in \autoref{fig:map},
where $7$ and $8$ are the integers in $A$'s or $D$'s.
\end{example}

\begin{lem}\label{thm:peq}
Let $p$ be a prime and $n=pt+r$ with $0\leq r<p$. Let
$H=\{s_1,s_2,\dots,s_h;E_1^{m_1}, \dots,E_\ell^{m_\ell}\}$ be an
equivalence class of $\S_{n,p}/\sim'$ described above. Then the number
of all permutations in the collection is
\begin{equation}\label{eqn:coll}
\left|\tilde f_p^{-1}(H)\right|=\frac{(1+(p-1)!)^{h}(p!)^{t-h} r!}
{m_1!m_2!\cdots m_\ell!}.
\end{equation}
\end{lem}
\begin{proof}
We need to enumerate the set $\tilde f_p^{-1}(H)$. Each permutation
in the set has the special disjoint cycle decomposition prescribed by $H$.
Recall that each $s_i$ can represent either a $p$-cycle or a $1$-cycle
of multicity $p$. 
If $s_i$ represents a $p$-cycle, it contributes a factor $(p-1)!$ to the
total number of permutations to be counted;
if it represents a $1$-cycle with multicity $p$, it contributes
a factor $1$. So in total each $s_i$ contributes a factor $(p-1)!+1$,
which explains the factor $(1+(p-1)!)^{h}$ in (\ref{eqn:coll}).

Now recall that $p$ is prime and $E$'s are a $p$-cycle or $1$-cycle.
Each $j\in[t]\setminus\{s_1,\ldots,s_h\}$
appears exactly $p$ times in $E_1^{m_1},\dots,E_\ell^{m_\ell}$,
which will be replaced by $p$ integers $p(j-1)+1,p(j-1)+2,\ldots,pj$,
contributing the factor $(p!)^{t-h}$ 
in (\ref{eqn:coll}); and $t+1$ appears exactly $r$ times, which
correspond to $r$ integers $pt+1,pt+2,\dots,pt+r$,
contributing a factor $r!$.
This argument overcounts the set $\tilde f_p^{-1}(H)$, since $E_i$ appears
$m_i$ times and the argument respects ordering of the cycles,
while we are interested in
unordered cycle decompositions.
Moreover, since each $E_i$ is a $1$-cycle or a $p$-cycle of type $B$, there
is no other repetition arising from a cyclic rotation inside a cycle in
$E_i$'s. So we need exactly the factor $\frac{1}{m_1!m_2!\cdots m_\ell!}$
in (\ref{eqn:coll}) to count the unordered structures.
\end{proof}

Now we can prove \eqref{grady} combinatorially.

\begin{thm}
Let $p$ be a prime and $n$ a positive integer. Then
$$\ord_p(\tau_p(n))\geq\floor np-\floor n{p^2}.$$
\end{thm}
\begin{proof}
Note that
$$
\tau_p(n)=|\S_{n,p}|=\sum_{H}\left|\tilde f_p^{-1}(H)\right|,
$$ where $H$ runs through all distinct equivalence classes of
$\S_{n,p}/\sim'$, i.e., distinct images of~$\tilde f_p$.  Thus it
suffices to show that for any equivalence class of $\S_{n,p}/\sim'$,
we have $\ord_p\left(\left|\tilde f_p^{-1}(H)\right|\right)\geq\floor
np-\floor n{p^2}$.

Let $H=\{s_1,s_2, \dots,s_h;E_1^{m_1},\dots,E_\ell^{m_\ell}\}$ be an
equivalence class of $\S_{n,p}/\sim'$. By \autoref{thm:peq}, we have
$$\left|\tilde f_p^{-1}(H)\right|=\frac{(1+(p-1)!)^{h}(p!)^{t-h} r!}
{m_1!m_2!\cdots m_\ell!}.$$
Since $(p-1)!\equiv-1\mod p$, $\ord_p$ of
the numerator is at least $t$. Moreover, $m_i\leq p$ for all $i$, and
if $m_i=p$ then $E_i$ is a $p$-cycle, which implies that there are at
most $\floor{n}{p^2}$ $m_i$'s with $m_i=p$. Thus we get
$\ord_p\left(\left|\tilde f_p^{-1}(H)\right|\right)\geq\floor np-
\floor n{p^2}$.
\end{proof}

\section{The power of 2 in the number of involutions}\label{tn}
For $p=2$, $\S_{n,p}$ is in fact the set of involutions in $\S_n$,
which will be denoted by~$\I_n$. Recall that $\inv_n$ stands for the
number of involutions in $\S_n$, i.e., $|\I_n|$. We will compute
$\ord_2(\inv_n)$ exactly and look at $\odd_n$ the odd factor of
$\inv_n$, i.e.,
$$\odd_n=\frac{\inv_n}{2^{\ord_2(\inv_n)}}.$$

Let $n=2t+r$ with $0\leq r<2$. Recall the equivalence relation $\sim'$
on $\S_{n,2}$ in \autoref{sec:comb}. Each equivalence class of
$\S_{n,2}/\sim'$ is represented by
$$H=\{s_1,s_2,\dots,s_h;E_1^{m_1},\dots,E_\ell^{m_\ell}\},$$ where $E$'s
denote either a $2$-cycle, consisting of two distinct integers, of
multiplicity at most two or a $1$-cycle with multiplicity one.  The
equivalence class may be represented by a graph $G=(\mathcal
V,\mathcal E)$ with vertex set
$$\mathcal V=\{v_1,v_2,\dots,v_t\},\mbox{ if $n=2t$};
\{v_1,v_2,\dots,v_{t+1}\},\mbox{ if $n=2t+1$},$$
and edge set $\mathcal E=\{\{a,b\}:\mbox{$(a,b)=E_j$, for some $j$ and
$a\ne b$}\}$, regarded as a multiset, where the multiplicity of the edge
corresponding to $E_j$ is $m_j$. We can construct $H$ from $G$
if we know $n$. 

Let $\G_n$ be the set of all graphs with vertex set
$$\{v_1,v_2,\dots,v_t\},\mbox{ if $n=2t$};
\{v_1,v_2,\dots,v_{t+1}\},\mbox{ if $n=2t+1$},$$
satisfying the following conditions:
\begin{itemize}
\item there is no loop,
\item the degree of each vertex is at most two, and that of
$v_{t+1}$ is at most one,
\item the multiplicity of each edge is at most two.
\end{itemize}
Then there is a one-to-one correspondence between the set of
equivalence classes of $\S_{n,2}/\sim'$ and the set $\G_n$.  Thus we
have the induced surjection $\tilde f_2:\I_n\rightarrow\G_n$.

Each connected component of a graph in $\G_n$ is either a cycle of
length at least two or a path.

The corollary below follows immediately from \autoref{thm:peq},
since a $2$-cycle is an edge with multiplicity $2$ in this case.
\begin{cor}\label{thm:2cycle}
Let $G\in\G_n$ have $s$ $2$-cycles. Then
$$\left|\tilde f_2^{-1}(G)\right|=2^{\floor n2-s}.$$
\end{cor}

The maximum number of $2$-cycles in a graph $G\in\G_n$ is $\floor n4$,
which gives $\ord_2(t_n)\geq
\floor n2-\floor n{4}$. Since there may be many such $G$'s, we need to
do more to determine
$\ord_2(t_n)$ exactly. Let $g_n$ denote the number of $G\in\G_n$ without
$2$-cycles. It is easy to see that
$$g_{2n+1}=g_{2n} + n g_{2n-1}.$$
For $n\leq 3$, $g_{2n}$ is just the number of simple (labeled) 
graphs with $n$ vertices. Thus $g_0=g_2=1$, $g_4=2$ and $g_6=8$.
Using the above recurrence, we get $g_1=1$, $g_3=2$, $g_5=6$ and $g_7=26$.
For more values of $g_n$, see \autoref{tab:+1}.

Let $\akn{a}{b}{n}$ denote the following product.
$$\akn{a}{b}{n}=\prod_{i=0}^{n-1}(a+ib)$$ Note that $\akn{1}{2}{n}$ is
always odd, in fact, it is the product of the first $n$ odd integers.

\begin{thm}\label{thm:tn}
Let $n=4k+r$ with $0\leq r<4$. Then
$$\inv_n=2^{k+\fl{r/2}}\sum_{i=0}^{k}2^{i}
\binom{k}{i}\frac{\akn{1}{2}{k+\fl{r/2}}}{\akn{1}{2}{i+\fl{r/2}}}\,g_{4i+r}.
$$
\end{thm}
\begin{proof}
Since $\tilde f_2:\I_n\rightarrow\G_n$ is a surjection, we have
$$\inv_n=\sum_{G\in\G_n}\left|\tilde f_2^{-1}(G)\right|.$$
If $G\in\G_n$ has $i$ $2$-cycles, then by \autoref{thm:2cycle},
$\left|\tilde f_2^{-1}(G)\right|=2^{\fl{n/2}-i}$.
Since the number of such $G$ is
$\binom{\fl{n/2}}{2i}\akn{1}{2}{i}\,g_{n-4i}$, we get
\begin{eqnarray*}
\inv_n&=&\sum_{i=0}^{k}2^{\fl{n/2}-i}\binom{\fl{n/2}}{2i}\akn{1}{2}{i}
g_{n-4i}\\
&=&\sum_{i=0}^{k}2^{\fl{n/2}-k+i}\binom{\fl{n/2}}{2k-2i}\akn{1}{2}{k-i}
g_{n-4k+4i}\\
&=&\sum_{i=0}^{k}2^{k+\fl{r/2}+i}\binom{2k+\fl{r/2}}{2i+\fl{r/2}}\akn{1}{2}{k-i}g_{4i+r}\\
&=&2^{k+\fl{r/2}}\sum_{i=0}^{k}2^{i}
\binom{k}{i}\frac{\akn{1}{2}{k+\fl{r/2}}}{\akn{1}{2}{i+\fl{r/2}}}
g_{4i+r}.
\end{eqnarray*}
\end{proof}

Since $g_0=g_1=g_2=1$, $g_3=2$ and $g_7=26$, we have the following
theorem, where $\delta_{r,3}$ is $1$, if $r=3$; $0$, otherwise.
\begin{thm}\label{2pwr}
Let $n=4k+r$ with $0\leq r<4$. Then the largest power of $2$ and the odd
factor $\odd_n$ of $\inv_n$ are the following:
\begin{align*}
\ord_2(\inv_n)&=k+\floor r2+\delta_{r,3}=\floor n2-2\floor n4+
\floor{n+1}4,\\
\odd_n&=\sum_{i=0}^{k}2^{i-\delta_{r,3}}\binom{k}{i}
\frac{\akn{1}{2}{k+\fl{r/2}}}{\akn{1}{2}{i+\fl{r/2}}}\,g_{4i+r}.
\end{align*}
\end{thm}

\section{Weighted sum of involutions}\label{sec:wt}
For $\pi\in\I_n$, let $\ai(\pi)$ denote the number of $i$-cycles in
$\pi$. We define the weight of an involution $\pi$ to be
$$\wt(\pi)=x^{\ao(\pi)} y^{\at(\pi)}.$$

Consider the weight generating function
\begin{equation}\label{eq:tnxy}
\inv_n(x,y)=\sum_{\pi\in\I_n}\wt(\pi).
\end{equation}
We can easily verify
$$\inv_n(x,y)= x\cdot\inv_{n-1}(x,y) + (n-1)y\cdot\inv_{n-2}(x,y).$$

Note that $\inv_n(x,-1)$ is the matchings polynomial of the complete
graph with $n$ vertices, which is equivalent to a Hermite
polynomial, see \cite{Godsil1993}.

We will find a formula for $\inv_n(x,y)$. Recall that $n=2t+r$ with
$0\leq r<2$ and the vertex set of a graph in $\G_n$ is either $[t]$
or $[t+1]$ depending on the parity of $n$. 
For $G\in\G_n$, we put the weight on each edge and vertex as follows:
\begin{itemize}
\item For every edge $e$, $\wt(e)=y$.
\item For $i\ne t+1$, $\wt(v_i)=\left\{
\begin{array}{ll}
1,&\mbox{if $\deg(v_i)=2$,}\\
x,&\mbox{if $\deg(v_i)=1$,}\\
\frac{x^2+y}2, &\mbox{if $\deg(v_i)=0$.}
\end{array}\right.$
\item $\wt(v_{t+1})=\left\{
\begin{array}{ll}
1,&\mbox{if $\deg(v_{t+1})=1$,}\\
x,&\mbox{if $\deg(v_{t+1})=0$.}
\end{array}\right.$
\end{itemize}
The weight $\wt(G)$ of $G$ is defined to be the product of weights of all
vertices and edges. It is not difficult to see that $\wt(G)$ is the
average of the weights of $\pi$ with $\tilde f_2(\pi)=G$, i.e.,
$$\sum_{\pi\in\tilde f_2^{-1}(G)}\wt(\pi)=|\tilde f_2^{-1}(G)|\wt(G).$$
Let
$$g_n(x,y) =\sum_G\wt(G),$$
where the sum is over all $G\in\G_n$ without 2-cycles.

Using the same argument in the proof of \autoref{thm:tn}, we have the
following theorem, since a $2$-cycle has two edges of weight $y$.

\begin{thm}\label{t_n}
Let $n=4k+r$ with $0\leq r<4$. Then
$$\inv_n(x,y)=2^{k+\fl{r/2}}\sum_{i=0}^{k}2^{i}
\binom{k}{i}\frac{\akn{1}{2}{k+\fl{r/2}}}{\akn{1}{2}{i+\fl{r/2}}}\,
y^{2k-2i}\,g_{4i+r}(x,y).$$
\end{thm}

We now find a recursion for $g_{n}(x,y)$.
\begin{prop}\label{rec}
Let $g_k(x,y)=0$ for negative integers $k$ and $g_0(x,y)=1$. Then for each positive
integer $n$, the following hold:
\begin{equation}\label{eq:g2n+1}
g_{2n+1}(x,y)=x\cdot g_{2n}(x,y)+ny\cdot g_{2n-1}(x,y),
\end{equation}
\begin{eqnarray}\label{eq:g2n}
    g_{2n}(x,y)&=&\frac{x^2+y}2 g_{2n-2}(x,y)+(n-1)xy\cdot
    g_{2n-3}(x,y)\\\notag
    &&+2\binom{n-1}{2} y^2\cdot g_{2n-4}(x,y)
    +3\binom{n-1}{3} y^4\cdot g_{2n-8}(x,y).
\end{eqnarray}
\end{prop}
\begin{proof}
\def\H{\mathfrak{H}}
The first recurrence, \eqref{eq:g2n+1}, is easy. For \eqref{eq:g2n},
let $\H_{2n}$ be the set of $G\in\G_{2n}$ without $2$-cycles.

We divide $\H_{2n}$ into four sets as follows:
\begin{align*}
\H_{2n}^{(0)} &=\{G\in\H_{2n}:\deg(v_n)=0\},\\
\H_{2n}^{(1)} &=\{G\in\H_{2n}:\deg(v_n)=1\},\\
\H_{2n}^{(2)} &=\{G\in\H_{2n}:\mbox{$v_n$ is contained in a 4-cycle}\},\\
\H_{2n}^{(*)} &=\{G\in\H_{2n}:\deg(v_n)=2
\mbox{ and $v_n$ is not contained in a 4-cycle}\}.
\end{align*}
Then it is easy to see that the weighted sums of $G$ in $\H_{2n}^{(0)}$,
$\H_{2n}^{(1)}$ and $\H_{2n}^{(2)}$ are, respectively, the first, second
and fourth terms in the right hand side of \eqref{eq:g2n}.

Let $G$ be a graph in $\H_{2n}^{(*)}$ and $a$, $b$ be the vertices
adjacent to $v_n$ in $G$. Let $G'$ denote the graph obtained from $G$
by collapsing the three vertices $v_n$, $a$ and $b$ 
to a new vertex $w$ as shown in \autoref{fig:collapse}.
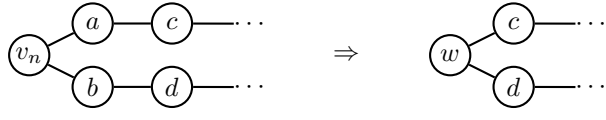
\begin{figure}
 \begin{center}
 \begin{pspicture}(0,0)(13,2)
\vput(0,1){n}
\xput(2,2){a}\xput(2,0){b}
\xput(4.5,2){c}\xput(4.5,0){d}
\rput(7,2){\rnode{A}{$\cdots$}}
\rput(7,0){\rnode{B}{$\cdots$}}
\ncline{n}{a}\ncline{n}{b}
\ncline{c}{a}\ncline{d}{b}
\ncline{c}{A}\ncline{d}{B}
\rput(10,1){$\Rightarrow$}
\end{pspicture}
\begin{pspicture}(0,0)(5,2)
\xput(0,1){w}
\xput(2,2){c}\xput(2,0){d}
\rput(4.5,2){\rnode{A}{$\cdots$}}
\rput(4.5,0){\rnode{B}{$\cdots$}}
\ncline{c}{w}\ncline{d}{w}
\ncline{c}{A}\ncline{d}{B}
\end{pspicture}
\caption{Collapsing $v_n$, $a$ and $b$ to $w$.}\label{fig:collapse}
 \end{center}
\end{figure}
Since $v_n$ is not contained in a $4$-cycle, there is no 2-cycle in
$G'$ and we can consider $G'$ as a graph in $\H_{2n-4}$ by relabeling
vertices.  Once $a,b$ and $w$ are fixed, for each $G'\in\H_{2n-4}$,
there are two graphs $G_1$ and $G_2$ in $\H_{2n}^{(*)}$ which collapse
to $G'$. For instance, if $w$ is an isolated vertex in $G'$, then $a$
and $b$ are connected to each other in $G_1$, and disconnected in
$G_2$. In this case, $\wt(G_1)=\wt(G')\frac{y^3}{(x^2+y)/2}$ and
$\wt(G_2)=\wt(G')\frac{y^2x^2}{(x^2+y)/2}$.
If $w$ is connected to $c$ and $d$ (one of them may be vacant),
then $a$ and $b$ are connected to $c$ and $d$ in $G_1$;
$d$ and $c$ in $G_2$ respectively. In this case, 
$\wt(G_1)=\wt(G_2)=y^2\wt(G')$. In both cases, we have
$\wt(G_1)+\wt(G_2)=2y^2\wt(G')$.
Thus the sum of $\wt(G)$ for $G\in\H_{2n}^{(*)}$ is equal to the third
term in the right hand side of \eqref{eq:g2n}.
\end{proof}
Using \autoref{rec}, we can compute $g_n(1,1)$ and $g_n(1,-1)$; see \autoref{tab:+1} and \autoref{tab:-1}. We will use these tables in the next section.

\begin{table}
 \begin{flushleft}
 \begin{tabular}{|c||c|c|c|c|c|c|c|c|c|c|c|c|c|c|c|}\hline
   $n$ & 0 & 1 & 2 & 3 & 4 & 5 & 6 & 7 & 8 & 9 & 10 & 11 & 12 &13 &14\\\hline
 $g_n$ & 1 & 1 & 1 & 2 & 2 & 6 & 8 & 26 & 41 & 145 & 253 & 978 & 1858 & 7726 & 15796
\\\hline
\end{tabular}
 \begin{tabular}{|c||c|c|c|c|c|c|c|c|c|c|c|c|c|c|c|}\hline
   $n$ & 15 & 16 & 17 & 18 & 19 &20 &21\\\hline
 $g_n$ &69878 &152219& 711243& 1638323& 8039510 & 99862594 & 252998224\\\hline
\end{tabular}
 \end{flushleft}
\caption{The values of $g_n=g_n(1,1)$ for $0\leq n\leq21$.}\label{tab:+1}
\end{table}

\begin{table}
\begin{flushleft}
 \begin{tabular}{|c||c|c|c|c|c|c|c|c|c|c|c|c|c|c|c|}\hline
   $n$ & 0 & 1 & 2 & 3 & 4 & 5 & 6 & 7 & 8 & 9 & 10 & 11 & 12 & 13 & 14 \\\hline
 $g_n(1,-1)$ & 1 & 1 & 0 & -1 & -1 & 1 & 2 & -1 & -6 & -2 & 
28 & 38 & -140 & -368 & 732 \\\hline
\end{tabular}
 \begin{tabular}{|c||c|c|c|c|c|c|c|}\hline
   $n$ & 15 & 16 & 17 & 18 & 19 & 20 & 21\\\hline
 $g_n(1,-1)$ & 3308 &-3934 & -30398 & 19232 & 292814
 & -44946 & -2973086\\\hline
%  -933740,
%  31770206,
%  25812556]
\end{tabular}
 \caption{The values of $g_n(1,-1)$ for $0\leq n\leq21$}\label{tab:-1}
 \end{flushleft}
\end{table}

\section{Odd and even involutions}\label{even}
Recall that $\sigma_2(\pi)$ is the number of $2$-cycles of $\pi$.
The {\it sign} of an involution $\pi\in\I_n$ is defined as usual, i.e.,
$$\s(\pi)=(-1)^{\at(\pi)}.$$ An involution is called {\it even}
(resp. {\it odd}), if the sign is 1 (resp. $-1$).  Let $\I_n^e$
(resp. $\I_n^o$) be the set of even (resp. odd) involutions in $\I_n$,
and let $\inv_n^e=|\I_n^e|$ and $\inv_n^o=|\I_n^o|$.

By definition of $\inv_n(x,y)$, we have
$$\inv_n(1,1)=\inv_n^e+\inv_n^o,\quad\inv_n(1,-1)=\inv_n^e-\inv_n^o.$$
Using the above equations, we will find $\ord_2(t_n^e)$ and $\ord_2(t_n^o)$. To do this we need
the following lemma.

\begin{lem}\label{two2}
Let $k$ and $i$ be positive integers. Then
$$\ord_2\left(2^i\binom{k}{i}\right)\geq\ord_2(k)+i-\ord_2(i).$$
Especially, we have
$$\ord_2\left(2^i\binom{k}{i}\right)\geq\ord_2(k)+1,$$
and if $i\geq5$, then
$$\ord_2\left(2^i\binom{k}{i}\right)\geq\ord_2(k)+3.$$
\end{lem}
\begin{proof}
It follows from the identity
$2^i\binom{k}{i}=2^i\cdot\frac{k}{i}\binom{k-1}{i-1}$.
\end{proof}

According to \autoref{t_n}, for $n=4k+r$ with $0\leq r<4$, we have
$$\inv_n(1,-1)=2^{k+\fl{r/2}}\sum_{i=0}^{k}2^{i}\binom{k}{i}
\frac{\akn{1}{2}{k+\fl{r/2}}}{\akn{1}{2}{i+\fl{r/2}}}\,g_{4i+r}(1,-1).$$
\begin{thm}
Let $n=4k+r$ with $0\leq r<4$. Then
$$\ord_2(\inv_n(1,-1))=\left\{
\begin{array}{ll}
k+\floor r2, &\mbox{if $r\ne 2$,}\\ k+3+\ord_2(k), &\mbox{if $r=2$.}
\end{array}\right.$$
\end{thm}
\begin{proof}
By \autoref{tab:-1}, we have $g_0(1,-1)=g_1(1,-1)=1, g_2(1,-1)=0$ and $g_3(1,-1)=-1$. Thus, if $r\ne 2$ then $\ord_2(\inv_n(1,-1))=\floor n2-\floor n4$.

If $r=2$, then $\inv_n(1,-1) = 2^{k+1}\sum_{i=0}^{k} a_i$ where
$a_i=2^{i}\binom{k}{i}\frac{\akn{1}{2}{k+1}}{\akn{1}{2}{i+1}} g_{4i+2}(1,-1)$.
Since $g_2(1,-1)=0$ and $g_6(1,-1)=2$, we have $a_0=0$ and
$\ord_2(a_1)=\ord_2(k)+2$. For $i\geq2$, using \autoref{tab:-1} and
\autoref{two2} we get $\ord_2(a_i)\geq\ord_2(k)+3$.  Thus
$\ord_2(\inv_{4k+2}(1,-1))=k+3+\ord_2(k)$.
\end{proof}
Now we can make a table of $\ord_2(\inv_n(1,1))$ and $\ord_2(\inv_n(1,-1))$; see \autoref{tab:main}.

\begin{table}
\begin{tabular}{|l|l|l|l|l|}\hline
$n$  & $\ord_2(\inv_n(1,1))$ & $\ord_2(\inv_n(1,-1))$ & $\ord_2(\inv_n^e)$ & $\ord_2(\inv_n^o)$\\\hline
$4k$ & $k$ & $k$ & $k+\chi_o(k)$& unknown\\
$4k+1$ & $k$ & $k$ & unknown & $k+\ord_2(k)+\chi_e(k)$\\
$4k+2$ & $k+1$ & $k+3+\ord_2(k)$ & $k$ & $k$\\
$4k+3$ & $k+2$ & $k+1$ & $k$&$k$\\
\hline
\end{tabular}
  \caption{The largest power of $2$ in the number of involutions, in the signed sum of involutions
and in the numbers of even or odd involutions.}\label{tab:main}
\end{table}

Since $\inv_n^e=\frac12(\inv_n(1,1)+\inv_n(1,-1))$ and
$\inv_n^o=\frac12(\inv_n(1,1)-\inv_n(1,-1))$, we get the following
corollary.
\begin{cor}
Let $k$ be a nonnegative integer. Then
$$\ord_2(\inv_{4k+2}^e) =\ord_2(\inv_{4k+2}^o)=
\ord_2(\inv_{4k+3}^e) =\ord_2(\inv_{4k+3}^o)=k.$$
\end{cor}

We find $\ord_2(\inv_{4k}^e)$ and $\ord_2(\inv_{4k+1}^o)$ in the
following two theorems separately. Let $\chi_o(n)$ (resp. $\chi_e(n)$) denote
1 if $n$ is odd (resp. even), and 0 otherwise.

\begin{thm}
Let $k$ be a nonnegative integer. Then
$$\ord_2(\inv_{4k}^e)=2\floor {k+1}2 = k+\chi_o(k).$$
\end{thm}
\begin{proof}
We have $\inv_{4k}^e = 2^{k}\sum_{i=0}^{k}a_i$, where
$$a_i= 2^{i-1}\binom{k}{i}\frac{\akn{1}{2}{k}}{\akn{1}{2}{i}}
\big(g_{4i}(1,1)+g_{4i}(1,-1)\big).$$
Using \autoref{tab:+1} and \autoref{tab:-1}, we have
\begin{align*}
g_{0}(1,1)+g_0(1,-1)&=1+1=2,\\
g_{4}(1,1)+g_4(1,-1)&=2-1=1,\\
g_{8}(1,1)+g_{8}(1,-1)&= 41-6\equiv3\mod4.
\end{align*}
Thus
$$ a_0 =\akn{1}{2}{k},\quad a_1 = k\akn{1}{2}{k},\quad a_2 = k(k-1)
\frac{\akn{1}{2}{k}}{3}\cdot(4q+3),$$ and
\begin{align*}
3(a_0+a_1+a_2) &=\akn{1}{2}{k}\left( 3+3k+(4q+3)(k^2-k)\right)\\
&\equiv\akn{1}{2}{k}\cdot 3(k^2+1)\mod 4.
\end{align*}
Thus $\ord_2(a_0+a_1+a_2)=\chi_o(k)$.  Since $\ord_2(a_i)\geq 2$ for
 $i\geq 3$, we finish the proof.
\end{proof}

\begin{thm}
Let $k$ be a nonnegative integer. Then
$$\ord_2(\inv_{4k+1}^o)=k+\ord_2(k)+\chi_e(k).$$
\end{thm}
\begin{proof}
We have $\inv_{4k+1}^o = 2^{k}\sum_{i=0}^{k}a_i$, where
$$a_i= 2^{i-1}\binom{k}{i}\frac{\akn{1}{2}{k}}{\akn{1}{2}{i}}
\big(g_{4i+1}(1,1) - g_{4i+1}(1,-1)\big).$$
Using \autoref{tab:+1} and \autoref{tab:-1}, we have
\begin{align*}
g_{1}(1,1)-g_{1}(1,-1) &= 1 -1 =0,\\
g_{5}(1,1)-g_{5}(1,-1) &= 6 -1 =5,\\
g_{9}(1,1)-g_{9}(1,-1) &= 145+2\equiv3\mod 4,\\
g_{17}(1,1)-g_{17}(1,-1) &= 711243 + 30398\equiv1\mod 2.
\end{align*}
Thus we can write $a_0=0$, $a_1=\akn{1}{2}{k}\cdot5k$,
$a_2=\akn{1}{2}{k}\binom{k}{2}\frac{2\cdot(4q_1+3)}3$,
$a_3=\akn{1}{2}{k}\binom{k}{3}\frac{2^2\cdot q_2}{5\cdot3}$ and
$a_4=\akn{1}{2}{k}\binom{k}{4}\frac{2^3\cdot(2q_3+1)}{7\cdot5\cdot3}$
for some integers $q_1,q_2$ and $q_3$.

Note that by \autoref{two2} we have
$\ord_2(a_i)\geq\ord_2(k)+3$ for $i\geq5$.
Thus, if $k$ is odd, then we have
$\ord_2(\inv_{4k+1}^o)=k$. 

Now assume that $k$ is even. Then
\begin{align*}
\ord_2(a_0+a_1+a_2) &=\ord_2(k(15+(k-1)(4q_1+3))),\\
\ord_2(a_3) &\geq\ord_2(k)+\ord_2(k-2)+1\geq\ord_2(k)+2,\\
\ord_2(a_4) &=\ord_2(k)+\ord_2(k-2).
\end{align*}
If $k=4m$, then $\ord_2(a_4)=\ord_2(k)+1$ and,
$\ord_2(a_0+a_1+a_2)\geq\ord_2(k)+2$.  If $k=4m+2$, then
$\ord_2(a_4)\geq\ord_2(k)+2$, and 
$\ord_2(a_0+a_1+a_2)=\ord_2(k)+1$.
Thus, if $k$ is even, then we always
have $\ord_2(a_0+\cdots+a_4)=\ord_2(k)+1$.

In all cases we have $\ord_2(\inv_{4k+1}^o)=k+\chi_e(k)(\ord_2(k)+1)=
k+\ord_2(k)+\chi_e(k)$.
\end{proof}

Now we can fill all the entries in \autoref{tab:main} except
$\ord_2(t_{4k+1}^e)$ and $\ord_2(t_{4k}^o)$. 
Based on Maple experiments, we conjecture the following.

\begin{conj}
There is a 2-adic integer $\rho=\sum_{i\geq0} \rho_i 2^i$,
with $0\leq\rho_i\leq1$, satisfying
$$\ord_2(\inv_{4k+1}^e)=k+\chi_o(k)\cdot(\ord_2(k+\rho)+1).$$
\end{conj}
For example, $\rho=1+2+2^3+2^8+2^{10}+\cdots$ satisfies the
condition for all $k\leq 1000$.

\section{The smallest period of $\odd_n\mod 2^s$}\label{period}
Chowla et al.~\cite{Chowla1951} proved that, if $m$ is odd, then
$\inv_{n+m}\equiv\inv_n\mod m$. We give their proof here for self
containment.

\begin{thm}\cite{Chowla1951}\label{thm:chowla}
If $m$ is odd, then $$\inv_{n+m}\equiv\inv_n\mod m.$$
\end{thm}
\begin{proof}
Induction on $n\geq0$.  We have
$$\inv_m=\sum_{2i+j=m}\frac{m!}{2^i i! j!} =
\sum_{2i+j=m}\frac{m!}{2^i (i+j)!}\binom{i+j}{j}\equiv 1\mod m,$$
because $\frac{m!}{2^i (i+j)!}\binom{i+j}{j}$ is divisible by $m$ if
$i>0$; and 1 if $i=0$.  Thus $\inv_{m+1}=\inv_m + m\inv_{m-1}\equiv 1
\mod m$.  We get $\inv_{n+m}\equiv\inv_n\mod m$ for $n=0,1$.
Suppose it holds for $n=0,1,\ldots,k$. Then it is true for $n=k+1$
because
\begin{align*}
\inv_{k+1+m}&=\inv_{k+m}+(k+m)\inv_{k+m-1}\\
&\equiv\inv_{k}+ k\inv_{k-1}\mod m\\
&=\inv_{k+1}.
\end{align*}
\end{proof}

The above theorem means that the sequence $\{\inv_n\mod m\}_{n\geq0}$
has a period $m$. In fact, $m$ is the smallest period.

\begin{thm}\label{thm:small}
Let $m$ be an odd integer. Then $m$ is the smallest period of the
sequence $\{\inv_n\mod m\}_{n\geq0}$.
\end{thm}
\begin{proof}
Let $d$ be the smallest period. Then
$\inv_{d}\equiv\inv_0\equiv1\mod m$,
$\inv_{d+1}\equiv\inv_1\equiv1\mod m$, and 
$\inv_{d+2}\equiv\inv_{2}\equiv 2\mod m$. On the other hand, we have
$\inv_{d+2}=\inv_{d+1}+(d+1)\inv_{d}\equiv d+2\mod m$.
Thus $m$ divides $d$, and we get $m=d$.
\end{proof}

If $m$ is even, then $\{\inv_{n}\mod m\}_{n\geq0}$ does not have
a period because $\inv_0=1$ but $\inv_{n}$ is even for all $n\geq2$.
However, there exists an integer
$N$ such that $\{\inv_{n}\mod m\}_{n\geq N}$ has a period. 

\begin{thm}
Let $\ell$ be an odd integer and $k$ be a positive integer.  Let
$m=2^k\ell$ and let $N$ be the smallest integer such that
$\{\inv_{n}\mod m\}_{n\geq N}$ has a period.  Then $N=4k-2$ and $\ell$
is the smallest period of $\{\inv_{n}\mod m\}_{n\geq N}$.
\end{thm}
\begin{proof}
By \autoref{2pwr}, we have $\ord_2(\inv_{4k-3})=k-1$ and
$\ord_2(\inv_n)\geq k$ for $n\geq 4k-2$.  Thus
$\inv_{4k-3+y}\not\equiv\inv_{4k-3}\mod 2^k$ for any positive integer $y$,
which implies $N\geq 4k-2$. On the other hand, we have
$\inv_{n+\ell}\equiv\inv_n\mod 2^k$ for $n\geq 4k-2$. 
Since $\inv_{n+\ell}\equiv\inv_n\mod \ell$ by
\autoref{thm:chowla}, we get $\inv_{n+\ell}\equiv\inv_n\mod m$ for
$n\geq 4k-2$. Thus $\{\inv_{n}\mod m\}_{n\geq 4k-2}$ has a period
$\ell$ and we get $N=4k-2$. 

It remains to show that $\ell$ is the smallest period. It is easy to see that
any period of $\{\inv_{n}\mod m\}_{n\geq N}$ is divisible by the smallest period 
of $\{\inv_{n}\mod \ell\}_{n\geq 0}$, which is $\ell$. Thus we get the theorem.
\end{proof}

Recall that $\odd_n$ is the odd factor of $t_n$. Similarly we can find
the smallest period of $\{\odd_n\mod 2^s\}_{n\geq0}$.  Let
$h(n)=\ord_2(t_n)=\floor{n}{2}-2\floor{n}{4}+\floor{n+1}4$.  Then
$\inv_n=2^{h(n)}\odd_{n}$. Thus we have
$$
\odd_{n+1}=2^{h(n)-h(n+1)}\odd_n+2^{h(n-1)-h(n+1)}n\odd_{n-1},
$$
which is equivalent to the following: if $n=4k+r$ with $0\leq r\leq3$ then
\begin{equation}\label{eq:odd}
\odd_{n+1}=2^{h(r)-h(r+1)}\odd_n+2^{h(r-1)-h(r+1)}n\odd_{n-1}.
\end{equation}

To find the smallest period of $\{\odd_n\mod 2^s\}_{n\geq0}$, we need the following two lemmas.
\begin{lem}\label{lem:oddfac}
Let $s\geq3$ be an integer. Then
$$\akn{1}{2}{2^{s-1}}\equiv1\mod 2^s.$$  
\end{lem}
\begin{proof}
Induction on $s$. It is true for $s=3$. Assume it is true for $s\geq3$.
Then $\akn{1}{2}{2^{s-1}}=2^sk+1$ for some integer $k$.
Then it holds for $s+1$ because
\begin{align*}
\akn{1}{2}{2^{s}}&=1\cdot3\cdot5\cdots(2^{s+1}-1)\\
&=\big(1\cdot3\cdot5\cdots(2^{s}-1)\big)\cdot
\big((2^{s+1}-1)(2^{s+1}-3)\cdots(2^{s+1}-(2^{s}-1))\big)\\
&\equiv\akn{1}{2}{2^{s-1}}\cdot(-1)^{2^{s-1}}\akn{1}{2}{2^{s-1}}\mod 2^{s+1}\\
&=2^{2s}k^2+2^{s+1}k+1\\
&\equiv1\mod 2^{s+1}.
\end{align*}
\end{proof}

\begin{lem}\label{s1}If $s\geq3$ then
$$\odd_{n+2^{s+1}}\equiv\odd_n\mod2^s.$$
\end{lem}
\begin{proof}
We use induction on $n$. First we will show that
$\odd_{2^{s+1}+n}\equiv 1\mod 2^s$ for $n=0,1$. By \autoref{2pwr},
\begin{align*}
\odd_{2^{s+1}+n}&=\sum_{i=0}^{2^{s-1}}2^{i}
\binom{2^{s-1}}{i}
\frac{\akn{1}{2}{2^{s-1}+\fl{n/2}}}{\akn{1}{2}{i+\fl{n/2}}}
\cdot\frac{g_{4i+n}}{2^{\delta_{n,3}}}
=\sum_{i=0}^{2^{s-1}}2^{i}
\binom{2^{s-1}}{i}
\frac{\akn{1}{2}{2^{s-1}}}{\akn{1}{2}{i}}g_{4i+n}.
\end{align*}
By \autoref{two2} and \autoref{lem:oddfac}, we get
$\odd_{2^{s+1}+n}\equiv\akn{1}{2}{2^{s-1}}\equiv1\mod 2^s$.

We have shown that the theorem is true for $n=0,1$. Assume $n\geq1$
and the theorem is true for all nonnegative integers less than $n+1$.
Then it is also true for $n+1$ because
if $n=4k+r$ for $0\leq r\leq3$ then by \eqref{eq:odd} we get
\begin{align*}
\odd_{n+1+2^{s+1}}&=2^{h(r)-h(r+1)}\odd_{n+2^{s+1}}+2^{h(r-1)-h(r+1)}(n+2^{s+1})\odd_{n-1+2^{s+1}}\\
&\equiv 2^{h(r)-h(r+1)}\odd_n+2^{h(r-1)-h(r+1)}n\odd_{n-1}\mod 2^s\\
&=\odd_{n+1}.
\end{align*}
\end{proof}

Now we have the following theorem.

\begin{thm}
If $s\geq3$ then $2^{s+1}$ is the smallest period of the sequence
$\{\odd_n\mod 2^s\}_{n\geq0}$.
\end{thm}
\begin{proof}
By \autoref{s1}, $2^{s+1}$ is a period. Since the smallest
period divides every period, it has to be $2^k$ for some $k$. It is
sufficient to show that $2^s$ is not a period.

Assume that $2^s$ is a period. By the recurrence relation
\eqref{eq:odd}, we have
$$\odd_{2^s+2} =\frac12\odd_{2^s+1} +\frac{2^s+1}2\odd_{2^s},\quad
\odd_{2^s+1} =\odd_{2^s} +2^s\cdot 2\odd_{2^s-1}.$$
Thus
$$\odd_{2^s+2} = (1+2^{s-1})\odd_{2^s} + 2^s\odd_{2^s-1}.$$ Since
$2^s$ is a period, $\odd_{2^s}\equiv\odd_0 =1\mod 2^s$.  Then we
have $\odd_{2^s+2}\equiv 1+2^{s-1}\mod 2^s$, which is a contradiction
to $\odd_{2^s+2}\equiv\odd_2 =1\mod 2^s$.
\end{proof}

\section*{Acknowledgement}
We would like to thank Professor Christian Krattenthaler for informing
us of the Ochiai's paper. We also thank the anonymous referees for
their careful reading and helpful comments.

\bibliographystyle{plain}

\end{document}